\documentclass[final,1p,times]{elsarticle}

\usepackage{amssymb}
 \usepackage{amsthm}
\usepackage{amscd}
\usepackage{amsmath}
\usepackage{amsfonts}
\usepackage{amssymb}
\usepackage{graphicx}
\newtheorem{theorem}{Theorem}

\newtheorem{lemma}[theorem]{Lemma}

\usepackage{mathrsfs}
\usepackage{titletoc}

\newcommand{\la}{\Delta}

\newcommand{\ra}{\rightarrow}

\newcommand{\f}{\frac}

\newcommand{\be}{\begin{equation}}
\renewcommand{\ra}{\rightarrow}
\newcommand{\ee}{\end{equation}}
\newcommand{\bea}{\begin{eqnarray}}
\newcommand{\eea}{\end{eqnarray}}
\newcommand{\bna}{\begin{eqnarray*}}
\newcommand{\ena}{\end{eqnarray*}}

\renewcommand{\le}{\left}
\newcommand{\ri}{\right}

\journal{***}

\begin{document}

\begin{frontmatter}
\title{Kazdan-Warner equation on graph}

\author{Alexander Grigor'yan}
\ead{grigor@math.uni-bielefeld.de}
\address{Department of Mathematics,
University of Bielefeld, Bielefeld 33501, Germany}

\author{Yong Lin}
 \ead{linyong01@ruc.edu.cn}

\author{Yunyan Yang}
 \ead{yunyanyang@ruc.edu.cn}

\address{ Department of Mathematics,
Renmin University of China, Beijing 100872, P. R. China}

\begin{abstract}
Let $G=(V,E)$ be a finite graph  and $\Delta$ be the usual
graph Laplacian. Using the calculus of variations and a method of upper and lower solutions, 
we give various conditions such that
the Kazdan-Warner equation
$\Delta u=c-he^u$
has a solution on $V$,
where $c$ is a constant, and $h:V\ra\mathbb{R}$ is a function. We also consider 
similar equations involving higher order derivatives on graph.
Our results can be compared with the original manifold case of Kazdan-Warner
(Ann. Math., 1974).
\end{abstract}

\begin{keyword}
 Kazdan-Warner equation on graph\sep Elliptic equation on graph

\MSC[2010] 34B45; 35A15; 58E30

\end{keyword}

\end{frontmatter}

\titlecontents{section}[0mm]
                       {\vspace{.2\baselineskip}}
                       {\thecontentslabel~\hspace{.5em}}
                        {}
                        {\dotfill\contentspage[{\makebox[0pt][r]{\thecontentspage}}]}
\titlecontents{subsection}[3mm]
                       {\vspace{.2\baselineskip}}
                       {\thecontentslabel~\hspace{.5em}}
                        {}
                       {\dotfill\contentspage[{\makebox[0pt][r]{\thecontentspage}}]}

\setcounter{tocdepth}{2}


\section{Introduction}

A basic problem in Riemannian geometry is that of describing curvatures on a given manifold.
Suppose that $(\Sigma,g)$ is a 2-dimensional compact Riemannian manifold without boundary, and $K$ is the
Gaussian curvature on it. Let $\widetilde{g}=e^{2u}g$ be a metric conformal to $g$, where
$u\in C^\infty(\Sigma)$. To find a smooth function $\widetilde{K}$ as the Gaussian curvature of
$(\Sigma, \widetilde{g})$, one is led to solving the nonlinear elliptic equation
\be\label{K-Z}\Delta_g u=K-\widetilde
{K}e^{2u},\ee
where $\Delta_g$ denotes the Laplacian operator on $(\Sigma,g)$. Let $v$ be a solution to
$\Delta_g v=K-\overline{K}$, where
$$\overline{K}=\f{1}{{\rm vol}_g(\Sigma)}\int_\Sigma Kdv_g.$$
Set $\psi=2(u-v)$. Then $\psi$ satisfies
$$\Delta \psi=2\overline{K}-(2\widetilde{K}e^{2v})e^\psi.$$
If one frees this equation from the geometric situation, then it is a special case of
\be\label{equation}\Delta_g u=c-he^u,\ee
where $c$ is a constant, and $h$ is some prescribed function, with neither $c$ nor $h$
depends on geometry of $(\Sigma,g)$. Clearly one can consider (\ref{equation}) in any dimensional
manifold. Now let $(\Sigma,g)$ be a compact Riemannian manifold of any dimension.
Note that the solvability of (\ref{equation}) depends on the sign of $c$. Let us summarize results of
Kazdan-Warner \cite{Kazdan-Warner}. For this purpose, think of
$(\Sigma,g)$ and $h\in C^\infty(\Sigma)$ as being fixed with ${\rm dim}\, \Sigma\geq 1$.\\

{\it Case 1. $c<0$.} A necessary condition for a solution is that $\overline{h}<0$, in which case there is a critical strictly
negative constant $c_ -(h)$ such that (\ref{equation}) is solvable if $c_-(h)<c<0$, but not solvable if $c<c_-(h)$.

{\it Case 2. $c=0$.} When ${\rm dim}\, \Sigma\leq 2$, the equation (\ref{equation}) has a solution if and only if both
$\overline{h}<0$ and $h$ is positive somewhere. When ${\rm dim}\, \Sigma\geq 3$, the necessary condition still holds.

{\it Case 3. $c>0$.} When ${\rm dim}\, \Sigma=1$, so that $\Sigma=S^1$, then (\ref{equation}) has a solution if and only if
$h$ is positive somewhere. When ${\rm dim}\, \Sigma=2$, there is a constant $0<c_+(h)\leq +\infty$ such that
(\ref{equation}) has a solution if $h$ is positive somewhere and if $0<c<c_+(h)$.\\

There are tremendous work concerning the Kazdan-Warner problem, among those we refer the reader to
Chen-Li \cite{Chen-Li1,Chen-Li2},  Ding-Jost-Li-Wang \cite{DJLW1,DJLW2}, and the references therein. \\

In this paper, we consider the Kazdan-Warner equation on a finite graph. In our setting, we shall prove the following:
In {\it Case 1}, we have the same conclusion as the manifold case; In {\it Case 2},  the equation
(\ref{equation}) has a solution if and only if both
$\overline{h}<0$ and $h$ is positive somewhere; While in {\it Case 3},  the equation
(\ref{equation}) has a solution if and only if $h$ is positive somewhere. Following the lines of
Kazdan-Warner \cite{Kazdan-Warner}, for results of {\it Case 2} and
{\it Case 3}, we use the
variational method; for results of {\it Case 1}, we use the principle of upper-lower solutions.
It is remarkable that Sobolev spaces on a finite graph
are all pre-compact. This leads to a very strong conclusion in {\it Case 3} compared with the manifold case.\\

We organized this paper as follows: In Section \ref{sec1}, we introduce some notations on graphs and state our main results.
In Section \ref{sec2}, we give two important lemmas, namely, the Sobolev embedding and the Trudinger-Moser embedding.
In Sections \ref{sec3}-\ref{sec5}, we prove Theorems \ref{Theorem 1}-\ref{Theorem 4} respectively. In Section \ref{sec6},
we discuss related equations involving higher order derivatives.

\section{Settings and main results}\label{sec1}

Let $G=(V,E)$ be a finite graph, where $V$ denotes the vertex set and $E$ denotes the edge set.
For any edge $xy\in E$, we assume that its weight $w_{xy}>0$ and that $w_{xy}=w_{yx}$.
Let $\mu:V\ra \mathbb{R}^+$ be a finite measure. For any function $u:V\ra \mathbb{R}$, the $\mu$-Laplacian (or Laplacian for short)
of $u$ is defined by
\be\label{lap}\Delta u(x)=\f{1}{\mu(x)}\sum_{y\sim x}w_{xy}(u(y)-u(x)),\ee
where $y\sim x$ means $xy\in E$.
The associated gradient form reads
\be\label{grad-form}\Gamma(u,v)(x)=\f{1}{2\mu(x)}\sum_{y\sim x}w_{xy}(u(y)-u(x))(v(y)-v(x)).\ee
Write $\Gamma(u)=\Gamma(u,u)$. We denote the length of its gradient by
\be\label{grd}|\nabla u|(x)=\sqrt{\Gamma(u)(x)}=\le(\f{1}{2\mu(x)}\sum_{y\sim x}w_{xy}(u(y)-u(x))^2\ri)^{1/2}.\ee
For any function $g:V\ra\mathbb{R}$,  an integral of $g$ over $V$ is defined by
\be\label{int}\int_V gd\mu=\sum_{x\in V}\mu(x)g(x),\ee
and an integral average of $g$ is denoted by
$$\overline{g}=\f{1}{{\rm{Vol}(V)}}\int_Vgd\mu=\f{1}{{\rm{Vol}(V)}}\sum_{x\in V}\mu(x)g(x),$$
where ${\rm Vol}(V)=\sum_{x\in V}\mu(x)$ stands for the volume of $V$.

The Kazdan-Warner equation on graph reads
\be\label{KZ-equ}\Delta u=c-he^u\quad{\rm in}\quad V,\ee
where $\Delta$ is defined as in (\ref{lap}), $c\in\mathbb{R}$, and $h:V\ra\mathbb{R}$ is a function.
If $c=0$, then (\ref{KZ-equ}) is reduced to
\be\label{0equ}\Delta u=-he^u\quad{\rm in}\quad V.\ee
Our first result can be stated as following:

 \begin{theorem}\label{Theorem 1}
 Let $G=(V,E)$ be a finite graph, and $h (\not\equiv 0)$ be a function on $V$.
 Then the equation (\ref{0equ}) has a solution if and only if $h$
 changes sign and
  $\int_Vhd\mu<0$.
 \end{theorem}

 In cases $c>0$ and $c<0$,  we have the following:

 \begin{theorem}\label{Theorem 2}
 Let $G=(V,E)$ be a finite graph, $c$ be a positive constant, and $h:V\ra\mathbb{R}$ be a function.
 Then the equation (\ref{KZ-equ}) has a solution if and only if $h$ is positive somewhere.
 \end{theorem}

 \begin{theorem}\label{Theorem 3}
 Let $G=(V,E)$ be a finite graph, $c$ be a negative constant, and $h:V\ra\mathbb{R}$ be a function.\\
 $(i)$ If (\ref{KZ-equ}) has a solution, then $\overline{h}<0$.\\
 $(ii)$ If $\overline{h}<0$, then there exists a constant $-\infty\leq c_-(h)<0$ depending on $h$ such that
 (\ref{KZ-equ}) has a solution for any $c_-(h)<c<0$, but has no solution for any $c<c_-(h)$.
  \end{theorem}

  Concerning the constant $c_-(h)$ in Theorem \ref{Theorem 3}, we have the following:

  \begin{theorem}\label{Theorem 4}
  Let $G=(V,E)$ be a finite graph, $c$ be a negative constant, and $h:V\ra\mathbb{R}$ be a function.
  Suppose that $c_-(h)$ is given as in Theorem \ref{Theorem 3}. If
  $h(x)\leq 0$ for all $x\in V$, but $h\not\equiv 0$, then $c_-(h)=-\infty$.
  \end{theorem}

 \section{Preliminaries}\label{sec2}

 Define a Sobolev space and a norm on it by
 $$W^{1,2}(V)=\le\{u:V\ra\mathbb{R}: \int_V(|\nabla u|^2+u^2)d\mu<+\infty\ri\},$$
 and
 $$\|u\|_{W^{1,2}(V)}=\le(\int_V(|\nabla u|^2+u^2)d\mu\ri)^{1/2}$$
 respectively.
 If $V$ is a finite graph, then $W^{1,2}(V)$ is exactly the set of all functions on $V$, a finite dimensional linear space.
 This implies the following Sobolev embedding:

 \begin{lemma}\label{Sobolv}Let $G=(V,E)$ be a finite graph.
 The Sobolev space $W^{1,2}(V)$ is pre-compact. Namely, if $\{u_j\}$ is bounded in $W^{1,2}(V)$, then there exists some
 $u\in W^{1,2}(V)$ such that up to a subsequence, $u_j\ra u$ in $W^{1,2}(V)$.
 \end{lemma}

 Also we have the following Trudinger-Moser embedding:

 \begin{lemma}\label{T-M-lemma}
 Let $G=(V,E)$ be a finite graph. For any $\beta>1$, there exists a constant $C$ depending only on
 $\beta$ and $V$ such that for all functions $v$ with $\int_V|\nabla v|^2d\mu\leq 1$ and $\int_Vvd\mu=0$,
 there holds
 $$\int_Ve^{\beta v^2}d\mu\leq C.$$
 \end{lemma}

 {\it Proof.} Let $\beta>1$ be fixed. For any function $v$ satisfying $\int_V|\nabla v|^2d\mu\leq 1$ and $\int_Vvd\mu=0$,
 we have by the Poincare inequality
 $$\int_Vv^2d\mu\leq C_0\int_V|\nabla v|^2d\mu\leq C_0,$$
 where $C_0$ is some constant depending only on $V$. Denote $\mu_{\min}=\min_{x\in V}\mu(x)$. In view of (\ref{int}),
 the above inequality leads to
 $\|v\|_{L^\infty(V)}\leq C_0/\mu_{\min}$. Hence
 $$\int_Ve^{\beta v^2}d\mu\leq e^{\beta C_0^2/\mu_{\min}}{\rm Vol}(V).$$
 This gives the desired result. $\hfill\Box$
 \section{The case $c=0$ }\label{sec3}

 In the case $c=0$, our approach comes out from that of Kazdan-Warner \cite{Kazdan-Warner}.

 {\it Proof of Theorem \ref{Theorem 1}.}

 {\it Necessary condition.} If (\ref{0equ}) has a solution $u$, then $e^{-u}\Delta u=-h$. Integration by parts
 gives
 \bna
 -\int_Vhd\mu&=&\int_Ve^{-u}\Delta ud\mu\\
 &=&-\int_V\Gamma(e^{-u},u)d\mu\\
 &=&-\f{1}{2}\sum_{x\in V}\sum_{y\sim x}w_{xy}(e^{-u(y)}-e^{-u(x)})(u(y)-u(x))\\
 &>&0,
 \ena
 since $(e^{-u(y)}-e^{-u(x)})(u(y)-u(x))\leq 0$ for all $x,y\in V$ and $u$ is not a constant.

 {\it Sufficient condition.} We use the calculus of variations. Suppose that $h$ changes sign and
 \be\label{mean-h}\int_Vhd\mu<0.\ee
 Define a set
 \be\label{B1}\mathcal{B}_1=\le\{v\in W^{1,2}(V): \int_Vhe^vd\mu=0,\,\,\int_Vvd\mu=0\ri\}.\ee
 We {\it claim} that 
 \be\label{B-non}\mathcal{B}_1\not=\varnothing.\ee To see this, since $h$ changes sign and (\ref{mean-h}),
 we can assume $h(x_1)>0$ for some $x_1\in V$. Take a function $v_1$ satisfying $v_1(x_1)=\ell$ and $v_1(x)=0$ for
 all $x\not= x_1$. Hence
 \bna
 \int_Vhe^{v_1}d\mu&=&\sum_{x\in V}\mu(x)h(x)e^{v_1(x)}\\
 &=&\mu(x_1) h(x_1)e^\ell+\sum_{x\not=x_1}\mu(x)h(x)\\
 &=&(e^\ell-1)\mu(x_1) h(x_1)+\int_Vhd\mu\\
 &>&0
 \ena
 for sufficiently large $\ell$. Writing
 $\phi(t)=\int_Vhe^{tv_1}d\mu$, we have by the above inequality that $\phi(1)>0$. Obviously
 $\phi(0)=\int_Vhd\mu<0$. Thus there exists a constant $0<t_0<1$ such that $\phi(t_0)=0$.
 Let $v^\ast=t_0v_1-\f{1}{{\rm vol}(V)}\int_Vt_0v_1d\mu$, where ${\rm vol}(V)=\sum_{x\in V}\mu(x)$ stands for
 the volume of $V$. Then $v^\ast\in \mathcal{B}_1$. This concludes our claim (\ref{B-non}).

 We shall minimize the functional
 $J(v)=\int_V|\nabla v|^2d\mu$.
 Let
 $$a=\inf_{v\in \mathcal{B}_1}J(v).$$
 Take a sequence of functions $\{v_n\}\subset\mathcal{B}_1$ such that $J(v_n)\ra a$. Clearly
 $\int_V|\nabla v_n|^2d\mu$ is bounded and $\int_V{v}_nd\mu=0$. Hence $v_n$ is bounded in $W^{1,2}(V)$.
 Since $V$ is a finite graph, the Sobolev embedding (Lemma \ref{Sobolv}) implies that up to a subsequence, $v_n\ra v_\infty$ in
 $W^{1,2}(V)$. Hence $\int_V{v_\infty}d\mu=0$, $\int_Vhe^{v_\infty}d\mu=\lim_{n\ra\infty}\int_Vhe^{v_n}d\mu=0$, and
 thus $v_\infty\in\mathcal{B}_1$. Moreover
 $$\int_V|\nabla v_\infty|^2d\mu=\lim_{n\ra\infty}\int_V|\nabla v_n|^2d\mu=a.$$
 One can calculate the Euler-Lagrange equation of  $v_\infty$ as follows:
 \be\label{E-L1}\Delta v_\infty=-\f{\lambda}{2}he^{v_\infty}-\f{\gamma}{2},\ee
 where $\lambda$ and $\gamma$ are two constants. Indeed, for any $\phi\in W^{1,2}(V)$, there holds
  \bea\nonumber
  0&=&\le.\f{d}{dt}\ri|_{t=0}\le\{\int_V|\nabla(v_\infty+t\phi)|^2d\mu-\lambda\int_Vhe^{v_\infty+t\phi}d\mu-
  \gamma\int_V(v_\infty+t\phi)d\mu\ri\}\\\nonumber
   &=&2\int_V\Gamma(v_\infty,\phi)d\mu-\lambda\int_Vhe^{v_\infty}\phi d\mu-\gamma\int_V\phi d\mu\\
   \label{euler}&=&-2\int_V(\Delta v_\infty)\phi d\mu-\lambda\int_Vhe^{v_\infty}\phi d\mu-\gamma\int_V\phi d\mu,
  \eea
  which gives (\ref{E-L1}) immediately.
    Integrating the equation (\ref{E-L1}), we have $\gamma=0$.
 We {\it claim} that $\lambda\not=0$. For otherwise, we conclude from $\Delta v_\infty=0$ and $\int_Vv_\infty d\mu=0$
 that $v_\infty\equiv 0\not\in \mathcal{B}_1$. This is a contradiction. We further {\it claim} that $\lambda>0$. This 
 is true because $\int_Vhd\mu<0$ and
 $$0<\int_V e^{-v_\infty}\Delta v_\infty d\mu=-\f{\lambda}{2}\int_Vhd\mu.$$
 Thus we can write $\f{\lambda}{2}=e^{-\vartheta}$ for some constant $\vartheta$. Then $u=v_\infty+\vartheta$ is a desired solution of (\ref{0equ}).
 $\hfill\Box$

 \section{The case $c>0$}\label{sec4}
  
 {\it Proof of Theorem \ref{Theorem 2}.}

 {\it Necessary condition.} Suppose $c>0$ and $u$ is a solution to (\ref{KZ-equ}). Since $\int_V\Delta ud\mu=0$, we have
 $$\int_Vhe^ud\mu=c{\rm Vol}(V)>0.$$
 Hence $h$ must be positive somewhere on $V$.\\

 {\it Sufficient condition.} Suppose $h(x_0)>0$ for some $x_0\in V$. Define a set
 $$\mathcal{B}_2=\le\{v\in W^{1,2}(V): \int_Vhe^vd\mu=c{\rm Vol}(V)\ri\}.$$
 We {\it claim} that $\mathcal{B}_2\not=\varnothing$. To see this, we set
 $$u_\ell(x)=\le\{\begin{array}{lll}
 \ell,\quad x=x_0\\[1.2ex]
 0,\quad x\not=x_0.
 \end{array}\ri.$$
 It follows that
 $$\int_Vhe^{u_\ell}d\mu\ra+\infty\quad{\rm as}\quad \ell\ra+\infty.$$
  We also set $\widetilde{u}_\ell\equiv -\ell$, which leads to
 $$\int_Vhe^{\widetilde{u}_\ell}d\mu=e^{-\ell}\int_Vhd\mu\ra 0\quad{\rm as}\quad \ell\ra+\infty.$$
 Hence there exists a sufficiently large $\ell$ such that $\int_Vhe^{u_\ell}d\mu>c{\rm Vol}(V)$ and
 $\int_Vhe^{\widetilde{u}_\ell}d\mu<c{\rm Vol}(V)$. We define a function $\phi:\mathbb{R}\ra\mathbb{R}$ by
 $$\phi(t)=\int_Vhe^{tu_\ell+(1-t)\widetilde{u}_\ell}d\mu.$$
 Then $\phi(0)<c{\rm Vol}(V)<\phi(1)$, and thus there exists a $t_0\in (0,1)$ such that $\phi(t_0)=c{\rm Vol}(V)$.
 Hence $\mathcal{B}_2\not=\varnothing$ and our claim follows.
 We shall solve (\ref{KZ-equ}) by minimizing the functional
 $$J(u)=\f{1}{2}\int_V|\nabla u|^2d\mu+c\int_Vud\mu$$
 on $\mathcal{B}_2$. For this purpose, we write $u=v+\overline{u}$, so $\overline{v}=0$. Then
 for any $u\in\mathcal{B}_2$, we have  $$\int_Vhe^vd\mu=c{\rm Vol}(V)e^{-\overline{u}}>0,$$ and thus
 \be\label{ju}J(u)=\f{1}{2}\int_V|\nabla u|^2d\mu-c{\rm Vol}(V)\log\int_Vhe^vd\mu+c{\rm Vol}(V)\log(c{\rm Vol}(V)).\ee
 Let $\widetilde{v}=v/\|\nabla v\|_2$. Then $\int_V\widetilde{v}d\mu=0$ and $\|\nabla \widetilde{v}\|_2=1$.
 By the Poincare inequality, $\|\widetilde{v}\|_2\leq C_0$ for some constant $C_0$ depending only on $V$.
 By Lemma \ref{T-M-lemma}, for any $\beta>1$, one can find a constant $C$ depending only on $\beta$ and $V$ such that
 \be\label{T-M}\int_Ve^{\beta {\widetilde{v}}^2}d\mu\leq C(\beta,V).\ee
 This together with an elementary inequality $ab\leq \epsilon a^2+\f{b^2}{4\epsilon}$ implies that for any $\epsilon>0$,
 \bna
 \int_Ve^vd\mu&\leq&\int_Ve^{\epsilon\|\nabla v\|_2^2+\f{v^2}{4\epsilon\|\nabla v\|_2^2}}d\mu\\
 &=&e^{\epsilon\|\nabla v\|_2^2}\int_Ve^{\f{v^2}{4\epsilon\|\nabla v\|_2^2}}d\mu\\
 &\leq& Ce^{\epsilon\|\nabla v\|_2^2},
 \ena
 where $C$ is a positive constant depending only on $\epsilon$ and $V$. Hence
 $$\int_Vhe^vd\mu\leq C(\max_{x\in V}h(x))e^{\epsilon\|\nabla v\|_2^2}.$$
 In view of (\ref{ju}), the above inequality leads to
 $$J(u)\geq \f{1}{2}\int_V|\nabla u|^2d\mu-c{\rm Vol}(V)\epsilon\|\nabla v\|_2^2-C_1,$$
 where $C_1$ is some constant depending only on $\epsilon$ and $V$. Choosing $\epsilon=\f{1}{4c{\rm Vol}(V)}$,
 and noting that $\|\nabla v\|_2=\|\nabla u\|_2$, we obtain for all $u\in\mathcal{B}_2$,
 \be\label{lowerb}J(u)\geq \f{1}{4}\int_V|\nabla u|^2d\mu-C_1.\ee
 Therefore $J$ has a lower bound on the set $\mathcal{B}_2$. This permits us to consider
 $$b=\inf_{u\in\mathcal{B}_2}J(u).$$
 Take a sequence of functions $\{u_k\}\subset\mathcal{B}_2$ such that $J(u_k)\ra b$. Let
 $u_k=v_k+\overline{u_k}$. Then $\overline{v_k}=0$, and it follows from (\ref{lowerb}) that
 $v_k$ is bounded in $W^{1,2}(V)$. This together with the equality
 $$\int_Vu_kd\mu=\f{1}{c}J(u_k)-\f{1}{2c}\int_V|\nabla v_k|^2d\mu$$
 implies that $\{\overline{u_k}\}$ is a bounded sequence. Hence $\{u_k\}$ is also bounded in $W^{1,2}(V)$.
 By the Sobolev embedding (Lemma \ref{Sobolv}), up to a subsequence,
 $u_k\ra u$ in $W^{1,2}(V)$. It is easy to see that $u\in\mathcal{B}_2$ and $J(u)=b$. Using
 the same method of (\ref{euler}), we derive the Euler-Lagrange equation of the minimizer $u$, namely,
 $\Delta u=c-\lambda he^u$
 for some constant $\lambda$. Noting that $\int_V\Delta ud\mu=0$, we have $\lambda=1$. Hence $u$ is a solution of
 the equation (\ref{KZ-equ}). $\hfill\Box$

\section{The case $c<0$}\label{sec5}

In this section, we prove Theorem \ref{Theorem 3} by using a method of upper and lower solutions. In particular,
we show that it suffices to construct an upper solution of the equation (\ref{KZ-equ}). This is exactly the graph version
of the argument of Kazdan-Warner (\cite{Kazdan-Warner}, Sections $9$ and $10$).\\

We call a function $u_-$  a lower solution of (\ref{KZ-equ}) if for all $x\in V$, there holds
$$\Delta u_-(x)-c+he^{u_-(x)}\geq 0.$$
Similarly, $u_+$ is called an upper solution of (\ref{KZ-equ}) if for all $x\in V$, it satisfies
$$\Delta u_+(x)-c+he^{u_+(x)}\leq 0.$$ We begin with the following:

\begin{lemma}\label{lemma9.3}
Let $c<0$. If there exist lower and upper solutions, $u_-$ and $u_+$, of the equation (\ref{KZ-equ}) with $u_-\leq u_+$,
then there exists a solution $u$ of (\ref{KZ-equ}) satisfying $u_-\leq u\leq u_+$.
\end{lemma}

{\it Proof.} We follow the lines of Kazdan-Warner (\cite{Kazdan-Warner}, Lemma 9.3). Set $k_1(x)=\max\{1,-h(x)\}$, so that
$k_1\geq 1$ and $k_1\geq -h$. Let $k(x)=k_1(x)e^{u_+(x)}$. We define $L\varphi\equiv\Delta\varphi-k\varphi$ and
$f(x,u)\equiv c-he^{u}$. Since $G=(V,E)$ is a finite graph and $\inf_{x\in V}k(x)>0$, we have that $L$ is a compact operator 
and ${\rm Ker}(L)=\{0\}$. Hence we can define inductively $u_{j+1}$
as the unique solution to
\be\label{equ1}Lu_{j+1}=f(x,u_j)-ku_j,\ee
where $u_0=u_+$. We {\it claim} that
\be\label{monotone}u_-\leq u_{j+1}\leq u_j\leq\cdots\leq u_+.\ee
To see this, we estimate
$$L(u_1-u_0)=f(x,u_0)-ku_0-\Delta u_0+ku_0\geq 0.$$
Suppose $u_1(x_0)-u_0(x_0)=\max_{x\in V}(u_1(x)-u_0(x))>0$. Then
$\Delta(u_1-u_0)(x_0)\leq 0$, and thus $L(u_1-u_0)(x_0)<0$. This is a contradiction. Hence $u_1\leq u_0$ on $V$.
Suppose $u_j\leq u_{j-1}$, we calculate
\bna
L(u_{j+1}-u_j)&=&k(u_{j-1}-u_j)+h(e^{u_{j-1}}-e^{u_j})\\
&\geq& k_1(x)(e^{u_+(x)}-e^{\xi})(u_{j-1}-u_j)\\
&\geq&0,
\ena
where $u_j\leq\xi\leq u_{j-1}$. Similarly as above, we have $u_{j+1}\leq u_j$ on $V$,
and by induction, $u_{j+1}\leq u_j\leq\cdots\leq u_+$ for any $j$. Noting that
$$L(u_--u_{j+1})\geq k(u_j-u_-)+h(e^{u_j}-e^{u_-}),$$
we also have by induction $u_-\leq u_j$ on $V$ for all $j$. Therefore (\ref{monotone}) holds. Since $V$ is finite, it is easy
to see that up to a subsequence,
$u_j\ra u$ uniformly on $V$. Passing to the limit $j\ra+\infty$ in the equation (\ref{equ1}),
one concludes that $u$ is a solution of (\ref{KZ-equ}) with
$u_-\leq u\leq u_+$. $\hfill\Box$\\

Next we show that the equation (\ref{KZ-equ}) has infinite lower solutions. This reduces the proof of Theorem \ref{Theorem 3} 
to finding its upper solution.

\begin{lemma}\label{lowersolution}
There exists a lower solution $u_-$ of (\ref{KZ-equ}) with $c<0$. Thus (\ref{KZ-equ}) has a solution if
and only if there exists an upper solution.
\end{lemma}

{\it Proof.} Let $u_-\equiv-A$ for some constant $A>0$. Since $V$ is finite, we have
$$\Delta u_-(x)-c+h(x)e^{u_-(x)}=-c+h(x)e^{-A}\ra -c\quad{\rm as}\quad A\ra+\infty,$$
uniformly with respect to $x\in V$.
Noting that $c<0$, we can find sufficiently large $A$ such that $u_-$ is a lower solution of (\ref{KZ-equ}).
$\hfill\Box$\\

{\it Proof of Theorem \ref{Theorem 3}}.

$(i)$ {\it Necessary condition.} If $u$ is a solution of (\ref{KZ-equ}), then
\bna
-\int_Vhd\mu&=&\int_Ve^{-u}\Delta ud\mu-c\int_Ve^{-u}d\mu\\
&=&-\int_V\Gamma(e^{-u},u)d\mu-c\int_Ve^{-u}d\mu\\
&>&0.
\ena

$(ii)$ {\it Sufficient condition.}
It follows from Lemmas \ref{lemma9.3} and \ref{lowersolution} that (\ref{KZ-equ}) has a solution
if and only if (\ref{KZ-equ}) has an upper solution $u_+$ satisfying
$$\Delta u_+\leq c-he^{u_+}.$$
Clearly, if $u_+$ is an upper solution for a given $c<0$, then $u_+$ is also an upper solution for all
$\widetilde{c}<0$ with $c\leq \widetilde{c}$. Therefore, there exists a constant $c_-(h)$ with $-\infty\leq
c_-(h)\leq 0$ such that (\ref{KZ-equ}) has a solution for any $c>c_-(h)$ but has no solution for
any $c<c_-(h)$.

We {\it claim} that $c_-(h)<0$ under the assumption $\int_Vhd\mu<0$. To see this, we let $v$ be a solution of
$\Delta v=\overline{h}-h$. There exists some constant $a>0$ such that
$$|e^{av}-1|\leq \f{-\overline{h}}{2\max_{x\in V} |h(x)|}.$$
Let $e^b=a$. If $c=\f{a\overline{h}}{2}$ and $u_+=av+b$, we have
\bna
\Delta u_+-c+he^{u_+}&=&ah(e^{av}-1)+\f{a\overline{h}}{2}\\
&\leq&a(\max_{x\in V}|h(x)|)|e^{av}-1|+\f{a\overline{h}}{2}\\
&\leq&\f{a\overline{h}}{2}-\f{a\overline{h}}{2}\\
&=&0.
\ena
Thus if $c={a\overline{h}}/{2}<0$, then the equation (\ref{KZ-equ}) has an upper solution $u_+$. Therefore, $\overline{h}<0$
implies that $c_-(h)\leq {a\overline{h}}/{2}<0$. $\hfill\Box$\\

{\it Proof of Theorem \ref{Theorem 4}.} We shall show that if $h(x)\leq 0$ for all $x\in V$, but $h\not\equiv 0$, then
(\ref{KZ-equ}) is solvable for all $c<0$. For this purpose, we let $v$ be a solution of $\Delta v=\overline{h}-h$.
Note that $\overline{h}<0$. Pick constants $a$ and $b$ such that $a\overline{h}<c$ and $e^{av+b}-a>0$. Let
$u_+=av+b$. Since ${h}\leq 0$,
\bna
\Delta u_+-c+he^{u_+}&=&a\Delta v-c+he^{av+b}\\
&=&a\overline{h}-ah-c+he^{av+b}\\
&\leq&h(e^{av+b}-a)\\
&\leq& 0.
\ena
Hence $u_+$ is an upper solution. Consequently, $c_-(h)=-\infty$ if $h\leq 0$ but $h\not\equiv 0$. $\hfill\Box$

\section{Some extensions}\label{sec6}

The equation (\ref{equation}) involving higher order differential operators was also extensively studied on manifolds,
 see for examples \cite{Djadli-Malchiodi,LLL} and the references therein.
 In this section, we shall extend Theorems \ref{Theorem 1}-\ref{Theorem 4} to nonlinear elliptic equations involving
higher order derivatives. For this purpose, we define the length of $m$-order gradient of $u$ by
 \be\label{high-deriv}|\nabla^mu|=\le\{\begin{array}{lll}
 |\nabla\Delta^{\f{m-1}{2}}u|,\,\,{\rm when}\,\,\, m
 \,\,{\rm is\,\,odd}\\[1.5ex]
 |\Delta^{\f{m}{2}}u|,\,\,{\rm when}\,\,\, m
 \,\,{\rm is\,\,even},
 \end{array}\ri.\ee
 where $|\nabla\Delta^{\f{m-1}{2}}u|$ is defined as in (\ref{grd}) for the function $\Delta^{\f{m-1}{2}}u$, and
 $|\Delta^{\f{m}{2}}u|$ denotes the usual absolute of the function $\Delta^{\f{m}{2}}u$. Define a Sobolev space by
 $$W^{m,2}(V)=\le\{v:V\ra\mathbb{R}: \int_V(|v|^2+|\nabla^mv|^2)d\mu<+\infty\ri\}$$
 and a norm on it by
 $$\|v\|_{W^{m,2}(V)}=\le(\int_V(|v|^2+|\nabla^mv|^2)d\mu\ri)^{1/2}.$$
 Clearly $W^{m,2}(V)$ is the set of all functions on $V$ since $V$ is finite. Moreover, we have
 the following Sobolev embedding and the Trudinger-Moser embedding:
 \begin{lemma}\label{m-Sob}
 Let $G=(V,E)$ be a finite graph. Then for any integer $m>0$, $W^{m,2}(V)$ is pre-compact.
 \end{lemma}
 \begin{lemma}\label{m-T-M-lemma}
 Let $G=(V,E)$ be a finite graph. Let $m$ be a positive integer. Then for any $\beta>1$, there exists a constant $C$ depending only on
 $m$, $\beta$ and $V$ such that for all functions $v$ with $\int_V|\nabla^m v|^2d\mu\leq 1$ and $\int_Vvd\mu=0$,
 there holds
 $$\int_Ve^{\beta v^2}d\mu\leq C.$$
 \end{lemma}

We consider an analog of (\ref{KZ-equ}), namely
\be\label{m-equ}\Delta^m u=c-he^u\quad{\rm in}\quad V,\ee
where $m$ is a positive integer, $c$ is a constant, and $h: V\ra\mathbb{R}$ is a function. Obviously (\ref{m-equ})
is reduced to (\ref{KZ-equ}) when $m=1$. Firstly we have the following:

\begin{theorem}\label{Theorem 5}
 Let $G=(V,E)$ be a finite graph, $h(\not\equiv 0)$ be a function on $V$, and
 $m$ be a positive integer.
 If $c=0$, $h$ changes sign, and $\int_Vhd\mu<0$,
 then the equation (\ref{m-equ}) has a solution.
 \end{theorem}

 {\it Proof.} We give the outline of the proof. Denote
 $$\mathcal{B}_3=\le\{v\in W^{m,2}(V):\, \int_Vhe^vd\mu=0,\,\int_Vvd\mu=0\ri\}.$$
 In view of (\ref{B1}), we have that $\mathcal{B}_3=\mathcal{B}_1$, since $V$ is finite. 
 Hence $\mathcal{B}_3\not=\varnothing$. Now we minimize the functional
 $J(v)=\int_V|\nabla^mu|^2d\mu$ on $\mathcal{B}_3$. The remaining part is completely analogous to that of the proof of
 Theorem \ref{Theorem 1}, except for replacing Lemma \ref{Sobolv} by Lemma \ref{m-Sob}.
 We omit the details but leave it to interested readers. $\hfill\Box$\\

 Secondly, in the case $c>0$, the same conclusion as Theorem \ref{Theorem 2} still holds for the equation (\ref{m-equ}) with $m>1$. 
 Precisely we have the following:

  \begin{theorem}\label{Theorem 6}
 Let $G=(V,E)$ be a finite graph, $c$ be a positive constant, $h:V\ra\mathbb{R}$ be a function,
 and $m$ be a positive integer.
 Then the equation (\ref{m-equ}) has a solution if and only if $h$ is positive somewhere.
 \end{theorem}
 {\it Proof.} Repeating the arguments of the proof of Theorem \ref{Theorem 2} except for 
 replacing Lemmas \ref{Sobolv} and \ref{T-M-lemma} by Lemmas \ref{m-Sob} and $\ref{m-T-M-lemma}$
 respectively, we get the desired result. $\hfill\Box$\\

 Finally, concerning the case $c<0$, we obtain a result weaker than Theorem \ref{Theorem 3}.

 \begin{theorem}\label{Theorem 7}
 Let $G=(V,E)$ be a finite graph, $c$ be a negative constant, $m$ is a positive integer,
 and $h:V\ra\mathbb{R}$ be a function such that
 $h(x)<0$ for all $x\in V$. Then the equation (\ref{m-equ}) has a solution.
 \end{theorem}

 {\it Proof.} Since the maximum principle can not be available for equations involving poly-harmonic operators,
   we use the calculus of variations instead of the method of upper and lower solutions. Let $c<0$ be fixed. Consider the functional
   \be\label{funct}J(u)=\f{1}{2}\int_V|\nabla^mu|^2d\mu+c\int_Vud\mu.\ee
   Set
   $$\mathcal{B}_4=\le\{u\in W^{m,2}(V): \int_Vhe^ud\mu=c{\rm Vol}(V)\ri\}.$$
   Using the same method of proving (\ref{B-non}) in the proof of Theorem \ref{Theorem 2}, we have $\mathcal{B}_4\not=\varnothing$.

   We now prove that $J$ has a lower bound on $\mathcal{B}_4$. Let $u\in \mathcal{B}_4$. Write $u=v+\overline{u}$.
   Then $\overline{v}=0$ and $$\int_Vhe^{v}d\mu=e^{-\overline{u}}c{\rm Vol}(V),$$ which leads to
   $$\overline{u}=-\log\le(\f{1}{c{\rm Vol}(V)}\int_Vhe^vd\mu\ri).$$
   Hence
   \be\label{lo-b}J(u)=\f{1}{2}\int_V|\nabla^mu|^2d\mu-c{\rm Vol}(V)\log\le(\f{1}{c{\rm Vol}(V)}\int_Vhe^vd\mu\ri).\ee
   Since $c<0$ and $h(x)<0$ for all $x\in V$, we have $\max_{x\in V}h(x)<0$, and thus
   \be\label{dl}\f{h}{c{\rm Vol}(V)}\geq \delta=\f{\max_{x\in V}h(x)}{c{\rm Vol}(V)}>0.\ee
   Inserting (\ref{dl}) into (\ref{lo-b}), we have
   \be\label{l-3}J(u)\geq \f{1}{2}\int_V|\nabla^mu|^2d\mu-c{\rm Vol}(V)\log\delta-c{\rm Vol}(V)\log\int_Ve^vd\mu.\ee
   By the Jensen inequality,
   \be\label{jens}\f{1}{{\rm Vol}(V)}\int_Ve^vd\mu\geq e^{\overline{v}}=1.\ee
   Inserting (\ref{jens}) into (\ref{l-3}), we obtain
   \be\label{lll}J(u)\geq \f{1}{2}\int_V|\nabla^mu|^2d\mu-c{\rm Vol}(V)\log\delta-c{\rm Vol}(V)\log{\rm Vol}(V).\ee
   Therefore $J$ has a lower bound on $\mathcal{B}_4$. Set
   $$\tau=\inf_{v\in \mathcal{B}_4}J(v).$$ Take a sequence of functions $\{u_k\}\subset\mathcal{B}_4$ such that
   $J(u_k)\ra \tau$. We have by (\ref{lll}) that
   \be\label{gra-b}\int_V|\nabla^mu_k|^2d\mu\leq C\ee
   for some constant $C$ depending only on $c$, $\tau$, $V$ and $h$. By (\ref{funct}), we estimate
   \be\label{gr-2}\le|\int_Vu_kd\mu\ri|\leq \f{1}{|c|}|J(u_k)|+\f{1}{2c}\int_V|\nabla^mu_k|^2d\mu.\ee
   The Poincare inequality implies that there exists some constant $C$ depending only on $m$ and $V$ such that
   \be\label{poinc}\int_V|u_k-\overline{u_k}|^2d\mu\leq C\int_V|\nabla^mu_k|^2d\mu\ee
   Combining (\ref{gra-b}), (\ref{gr-2}), and (\ref{poinc}), one can see that $\{u_k\}$ is bounded in $W^{m,2}(V)$.
   Then it follows from Lemma \ref{m-Sob} that there exists some function $u$ such that up to a subsequence,
   $u_k\ra u$ in $W^{m,2}(V)$. Clearly $u\in\mathcal{B}_4$ and $J(u)=\lim_{k\ra\infty}J(u_k)=\tau$. In other words,
   $u$ is a minimizer of $J$ on the set $\mathcal{B}_4$. It is not difficult to check that (\ref{m-equ}) is the Euler-Lagrange equation 
   of $u$.
   This completes the proof of the theorem. $\hfill\Box$\\

 {\bf Acknowledgements.} A. Grigor'yan is partly supported  by SFB 701 of the German Research Council. Y. Lin is supported by the National Science Foundation of China (Grant No.11271011). Y. Yang is supported by the National Science Foundation of China (Grant No.11171347).


\bigskip

\begin{thebibliography}{00}
\bibitem{Chen-Li1} W. Chen, C. Li, Qualitative properties of solutions to some nonlinear elliptic
equations in $\mathbb{R}^2$, Duke Math. J. 71 (1993) 427-439.

\bibitem{Chen-Li2} W. Chen, C. Li, Gaussian curvature on singular surfaces, J. Geom. Anal. 3 (1993) 315-334.

\bibitem{DJLW1} W. Ding, J. Jost, J. Li, G. Wang,
       The differential equation $\la u=8\pi-8\pi he^u$ on a compact
        Riemann Surface, Asian J. Math. 1 (1997), 230-248.

\bibitem{DJLW2} W. Ding, J. Jost, J. Li, G. Wang,
       An analysis of the two-vortex case in the Chern-Siomons Higgs model, Calc. Var. 7 (1998), 87-97.
\bibitem{Kazdan-Warner}J. Kazdan, F. Warner, Curvature functions for compact $2$-manifolds,
   Ann. of Math. (2) 99 (1974), 14-47.

 \bibitem{Djadli-Malchiodi} Z.  Djadli, A. Malchiodi, Existence of conformal metrics with constant $Q$-curvature,
 Ann. Math. (2) 168  (2008) 813-858.

 \bibitem{LLL} J. Li, Y. Li, P. Liu, The $Q$-curvature on a 4-dimensional Riemannian manifold $(M,g)$  with $\int_M Qdv_g =8\pi^2$,
 Adv. Math.  231  (2012) 2194-2223.

\end{thebibliography}
\end{document}